\documentclass[12pt,a4paper]{article}
\newtheorem{theoreme}{Th\'eor\`eme}[section]

\newtheorem{lemme}[theoreme]{Lemme}
\newtheorem{corollaire}[theoreme]{Corollaire}

\newtheorem{remarque}[theoreme]{Remarque}

\newenvironment{demonstration}{\par\medskip{\it D\'emonstration }}{\quad\nobreak\vrule height 6pt width 6pt}
\makeatletter\@afterindenttrue\makeatother
\begin{document}

\title{Une infinit\'e de structures de contact
tendues sur les vari\'et\'es toro\"\i dales
\footnote{ Classification AMS : 57R17, 53D10.
Mots cl\'es : vari\'et\'e de contact, structure tendue, torsion, vari\'et\'e
toro\"\i dale.}}

\author{Vincent Colin}
\date{Octobre 2000}
\maketitle
\begin{abstract}
\vskip 0.3cm
We show that every closed toroidal irreducible orientable
$3$-manifold carries infinitely many universally tight contact
structures.

\centerline{\bf{R\'esum\'e}}
\vskip 0.3cm
On d\'emontre que toute vari\'et\'e orientable, irr\'eductible,
close et toro\"\i dale de dimension~$3$ porte une infinit\'e de structures
de contact universellement tendues.
\end{abstract}

\section{Introduction}
Une vari\'et\'e de dimension $3$ est dite
{\it toro\"\i dale} si elle contient un
tore {\it incompressible}, c'est-\`a-dire plong\'e et $\pi_1$-inject\'e.
L'objectif de ce texte est de d\'emontrer le
th\'eor\`eme suivant, qui r\'epond partiellement \`a une conjecture
d'E.~Giroux~\cite{[Gi5], [Gi6]}~:
\begin{theoreme}\label{theoreme} : Toute vari\'et\'e de dimension $3$,
orientable, irr\'eductible, close
(compacte sans bord) et toro\"\i dale porte
une infinit\'e de structures de contact universellement tendues deux \`a deux
non isomorphes.
\end{theoreme}

Jusqu'ici, le principal r\'esultat \'etabli
l'\'etait pour les vari\'et\'es irr\'eductibles
contenant un tore normal \cite{[Co3]}~:
un tore est {\it normal} s'il est
incompressible et si on ne peut pas
le disjoindre par isotopie d'un autre
tore incompressible. Cette propri\'et\'e
caract\'erise des vari\'et\'es dans lesquelles
se plongent de mani\`ere $\pi_1$-injective
un fibr\'e de Seifert de base ``assez large".
Dans \cite{[Co3]}, on d\'emontre ainsi
que toute vari\'et\'e orientable et irr\'eductible
qui contient un tore normal porte une infinit\'e
de structures de contact universellement tendues.

Comme dans ce premier travail,
la d\'emonstration du th\'eor\`eme~\ref{theoreme}
repose sur la notion de torsion, introduite
par E.~Giroux~\cite{[Gi3],[Gi4],[Gi5]}, et sur l'\'etude qu'il en a faite
dans le cas du tore \'epais.\\
\\
Avant de poursuivre plus avant la discussion, on rappelle quelques
notions de g\'eom\'etrie de contact.

Une {\it structure de contact} sur une vari\'et\'e $V$ de dimension $3$
est un champ de plans $\xi$ qui poss\`ede une \'equation
locale $\alpha =0$ telle que la $3$-forme $\alpha \wedge d\alpha$
soit une forme de volume. En particulier, sur une vari\'et\'e
orient\'ee, toute structure de contact poss\`ede un signe (donn\'e
par celui de $\alpha \wedge d\alpha$ qui se r\'ev\`ele
ind\'ependant du choix de $\alpha$).
Par d\'efaut, toutes les vari\'et\'es de contact rencontr\'ees
par la suite sont suppos\'ees orient\'ees
et les structures de contact positives.

Les structures de contact se
scindent en deux cat\'egories compl\'ementaires~:
les structures {\it vrill\'ees} et les structures {\it tendues}.
Une structure $\xi$ sur une vari\'et\'e
$V$ est vrill\'ee s'il existe un disque $D$ plong\'e
dans $V$ tel que $\xi \vert_{\partial D} =T_{\partial D} D$.
Si le rappel dans tout rev\^etement de $V$ d'une structure
vrill\'ee est une structure
vrill\'ee, ce n'est pas toujours le cas des structures
tendues (voir par exemple~\cite{[Co2]}) ; c'est
pourquoi on distingue comme sous-cat\'egorie de ces
derni\`eres les structures {\it universellement tendues}, qui
restent tendues apr\`es un passage au rev\^etement universel.
Aujourd'hui, les structures vrill\'ees sont parfaitement
comprises gr\^ace au travail de Y.~Eliashberg~\cite{[El1]}. Les
structures tendues en revanche restent en partie myst\'erieuses.
Leur \'etude est d'autant plus int\'eressante qu'elle se situe
\`a l'intersection de
nombreuses branches de la g\'eom\'etrie~:
th\'eorie des n\oe uds, th\'eorie des feuilletages,
g\'eom\'etrie symplectique...

On d\'ecrit \`a pr\'esent des outils intervenant
dans l'\'etude des vari\'et\'es de contact.

Si $S$ est une surface orient\'ee dans $V$ et si
$\xi$ est orient\'ee pr\`es de $S$,
en dehors des points $(x_i)_{i\in I}$ de  $S$,
g\'en\'eriquement isol\'es, o\`u $\xi (x_i ) =T_{x_i} S$,
l'intersection de $\xi$ et de $TS$ est un champ de droites
orient\'ees qui s'int\`egre en un feuilletage de $S$, appel\'e
{\it feuilletage caract\'eristique} de $S$ et not\'e $\xi S$.
Il est singulier en les points $x_i$.
Le feuilletage caract\'eristique trac\'e par une structure
$\xi$ sur une surface $S$ d\'etermine $\xi$ pr\`es de $S$.
En particulier, pour recoller deux structures de contact de
m\^eme
signe situ\'ees de part et d'autre d'une surface,
il suffit que leurs feuilletages caract\'eristiques co\"\i ncident.

Toute courbe int\'egrale de $\xi$ est dite {\it legendrienne}.
Si une courbe legendrienne orient\'ee $\gamma$ est contenue dans une surface
$S$, on d\'efinit son {\it invariant de Thurston-Bennequin relatif},
not\'e $tb(\gamma ,S)$, comme la moiti\'e de la somme alg\'ebrique, le
long de $\gamma$, des intersections entre $\xi$ et $TS$.
Si $\gamma$ est une courbe positivement transverse \`a $\xi$ et borde
une surface orient\'ee $S$, le champ
de plans $\xi \vert_S$ admet une section non singuli\`ere
$X$. On appelle {\it autoenlacement} de $\gamma$, not\'e $l(\gamma )$,
l'enlacement entre $\gamma$ et toute courbe obtenue
en poussant un peu $\gamma$ par $X$.

Pour diff\'erencier les structures de contact tendues,
E.~Giroux~\cite{[Gi3], [Gi4], [Gi5]} a introduit la notion de {\it torsion}.
Si $T\subset (V,\xi)$ est un tore incompressible
dans une vari\'et\'e de contact,
on cherche a plonger, pour le plus
grand entier naturel $n$ possible, le
produit de contact $(T^2 \times [0,2\pi ] ,\ker (\cos ntdx +\sin ntdy))$
dans la m\^eme classe d'isotopie $C$ que celle d'un voisinage
tubulaire de $T$.
Cet entier (\'eventuellement nul, si de tels plongements
n'existent pas, ou infini) est la torsion de $\xi$
dans la classe $C$ et est not\'e $Tor (V,\xi ,C)$,
ou $Tor (V,\xi )$ s'il n'y a pas d'ambigu\"\i t\'e.
On a vu dans \cite{[Co3]} comment lire cet
invariant sur le nombre de Thurston-Bennequin relatif
d'un certain type de courbe.
C'est cette strat\'egie que l'on va \`a nouveau appliquer
pour d\'emontrer le th\'eor\`eme~\ref{theoreme}.

Pour de plus amples d\'etails sur les notions de base de la g\'eom\'etrie de
contact, et notamment sur la cr\'eation et l'\'elimination
de singularit\'es d'un feuilletage $\xi S$ par isotopie de
$S$, on renvoie \`a \cite{[Gi2]}.

Je remercie chaleureusement E.~Giroux de l'aide
qu'il m'a apport\'ee dans la recherche de cette d\'emonstration.
Je lui dois notamment l'id\'ee de se placer sur
un rev\^etement de groupe $\pi_1 (T)$ ainsi que de distinguer
les diff\'erentes possibilit\'es pour le rappel de $T$.
Une partie de la r\'edaction de ce texte a \'et\'e effectu\'ee
lors d'un s\'ejour \`a l'universit\'e de Stanford
et \`a l'American Institute of Mathematics.
Je remercie ces deux institutions de leur soutien.

\section{Recueil de r\'esultats sur la torsion}
L'objet de cette section est de r\'epertorier
des r\'esultats, principalement dus \`a E.~Giroux, K.~Honda et Y.~Eliashberg, concernant
la torsion des structures de contact universellement tendues
sur le tore \'epais, l'anneau \'epaissi et le demi-espace.
On donne si n\'ecessaire des esquisses de preuves
permettant d'expliciter le lien \`a ces travaux
lorsqu'il n'est pas totalement transparent.
Le lemme~~\ref{lemme : fini} peut \^etre
\'egalement d\'eduit des arguments d\'evelopp\'es dans \cite{[Co3]}.

\subsection{Structures de contact universellement tendues
sur le tore \'epais}
{\bf Calcul de la torsion sur un mod\`ele}

\begin{lemme}\label{lemme : torsion}{\rm \cite{[Gi5]}} : Soit $\xi$ la structure
de contact d\'efinie sur $T^2 \times [-1,1]=\lbrace (x,y,\theta )\rbrace$
par l'\'equation $\cos f(\theta )dx+\sin f(\theta )dy=0$, o\`u
$f:[-1,1] \rightarrow R$ est \`a d\'eriv\'ee strictement positive
et v\'erifie~: $f(-1) \in [0,2\pi [$ et $f(1)\in [2n\pi ,2(n+1) \pi[$, $n\in N$.
Alors $Tor (T^2 \times [-1,1] ,\xi )$ est la partie enti\`ere de
$(f(1)-f(-1))/2\pi$, c'est-\`a-dire $n$ ou $n-1$.
\end{lemme}
$\;$\\
{\bf Addition des torsions}

\begin{lemme}\label{lemme : addition}{\rm \cite{[Gi5]}} :
Soit $\xi$ une structure de contact universellement
tendue sur $T^2 \times R$ telle que $\xi T^2 \times \lbrace \pm 1\rbrace$
soit un feuilletage lin\'eaire.
On a~: $$Tor (T^2 \times R,\xi \vert_{T^2 \times R} )=Tor(T^2 \times ] -\infty ,-1 ]
,\xi \vert_{T^2 \times ]-\infty ,-1]})
+Tor(T^2 \times [-1,1] ,\xi \vert_{T^2 \times [ -1,1]})$$$$+
Tor (T^2 \times [1,+\infty [ ,\xi \vert_{T^2 \times [1,+\infty [} )+k,$$
o\`u $k=0,1,$ ou $2$.
\end{lemme}
$\;$\\
{\bf Crit\`ere d'annulation de la torsion}

\begin{lemme}\label{lemme : 0} :
Soit $\xi$ une structure universellement
tendue sur $T^2 \times [-1,1]$ qui trace sur $T^2 \times \lbrace 1\rbrace$ un feuilletage
sans singularit\'e ni composante de Reeb
avec un nombre fini d'orbites p\'eriodiques non d\'eg\'en\'er\'ees
et un feuilletage lin\'eaire sur
$T^2 \times \lbrace -1\rbrace$.
On note $A=S^1 \times \lbrace * \rbrace \times [-1,1] \subset T^2 \times [-1,1]$.
On suppose que $\partial A$ est transverse,
ou \'eventuellement tangent le long de
$T^2 \times \lbrace -1\rbrace \cap \partial A$, \`a $\xi$.
S'il existe un arc legendrien $\gamma$ plong\'e dans $A$ joignant
un bord de $A$ \`a l'autre et tel que $ tb(\gamma ,A) =- \frac12$,
alors $Tor (T^2 \times [-1,1] ,\xi )=0$.
\end{lemme}
\begin{demonstration} : On peut toujours supposer
que, quitte \`a changer $A$ en conservant $\gamma$,
le bord de $A$ est parall\`ele \`a aux
orbites de $\xi T^2 \times \lbrace 1\rbrace$. Apr\`es cette modification, on est assur\'e
que $tb(\gamma, A)=0$, $-\frac12$ ou $-1$.
D\`es lors, pour tout feuilletage non singulier
$F$ de $T^2 \times \lbrace 1 \rbrace$, transverse \`a $\partial A$,
dont les orbites ferm\'ees sont non d\'eg\'en\'er\'ees
et co\"\i ncident avec celles de $A$ (avec une orientation
\'eventuellement diff\'erente, mais les orbites attractives restent attractives
et les r\'epulsives, r\'epulsives),
il existe une isotopie $C^0$-petite de
$T^2 \times \lbrace 1 \rbrace$ fixant $A$ (\cite{[Gi1]}),
telle que le feuilletage caract\'eristique du nouveau
tore
soit $F$ (ce feuilletage
comporte \'eventuellement des composantes de Reeb).
Ce faisant, d'apr\`es \cite{[Gi5]},
 la torsion du produit d\'elimit\'e par
$T^2 \times \lbrace -1\rbrace$ et l'image de $T^2 \times \lbrace 1\rbrace$
est la m\^eme que celle de $(T^2 \times [-1,1] ,\xi )$.
On note $(T^2 \times [-1,1] ,\xi_F )$ ce nouveau produit.
Par construction, on a toujours l'anneau $A$ contenant l'arc legendrien
$\gamma$ qui joignent tous deux un bord  de $(T^2 \times [-1,1] ,\xi_F )$
\`a l'autre.
On a de plus $tb(\gamma, A)=0$, $-\frac12$ ou $-1$.

Toujours d'apr\`es \cite{[Gi5]} et sa description
compl\`ete des structures de contact tendues sur le
tore \'epais, si $n$ d\'esigne la torsion de $\xi$,
il existe un feuilletage $F$ de $T^2 \times \lbrace 1 \rbrace$, un entier $k$,
un r\'eel $\theta_0$ et une fonction $f$ strictement
sup\'erieure \`a $2n\pi$
tels que $(T^2 \times [-1,1] ,\xi_F )\simeq \lbrace (x,y,\theta \leq f(y))\rbrace
\subset (T^2 \times [0,2(n+k)\pi ],\ker (\cos (\theta +\theta_0 )dx+\sin (\theta +\theta_0 )dy ))$.
Dans ces coordonn\'ees, on prolonge $A$ par l'anneau produit
$(\partial A\cap \lbrace \theta =f(y)\rbrace )\times \lbrace f(y) \leq \theta \leq
2(n+k)\pi \rbrace$, et l'arc $\gamma_F$ par
$(\partial \gamma_F \cap \lbrace \theta =f(y)\rbrace )\times \lbrace f(y) \leq \theta \leq
2(n+k)\pi \rbrace$.
On obtient ainsi un anneau $A_F$ portant un arc legendrien
$\gamma_F$ qui joignent tous deux un bord \`a l'autre
de $T^2 \times [0,2(n+k)\pi ]$.
Par construction, $tb(\gamma_F ,A_F) \geq -k-\frac12$
et, comme $\partial A_F$ est transverse \`a $\xi$,
on peut \'egalement supposer, quitte
\`a effectuer une isotopie de $A_F$ relative
\`a $\gamma_F$ ne modifiant pas $tb(\gamma_F ,A_F)$, que
$\partial A_F$ est transverse
(ou \'eventuellement partout tangent) \`a la structure.

D'apr\`es \cite{[Ka]} (th\'eor\`eme~7.6), on a alors $tb(\gamma_F ,A_F )\leq -n-k$.
On en d\'eduit que $n\leq \frac12$, ce qui implique que l'entier naturel $n$
est nul.
\end{demonstration}

$\;$\\
{\bf Lemme de r\'ealisation}

\begin{lemme}\label{lemme : normalisation}{\rm \cite{[Gi5]}} :
Soit $\xi$ une structure de contact universellement tendue
sur $T^2 \times [0,+\infty [$ qui trace un feuilletage lin\'eaire
sur $T^2 \times \lbrace 0\rbrace$.
Si $\xi$ est de torsion finie $n$, il existe $(c,\theta_0 )\in R^2$ tel que
pour tout $k\in N$, il existe un plongement
$\phi$ de $T^2 \times [0,k ]$ dans $T^2 \times [0,c]$
laissant invariant $T^2 \times \lbrace 0\rbrace$, et tel que $\phi_* \xi$
ait pour \'equation $\cos (\theta +\theta_0 ) dx+\sin (\theta +\theta_0 )dy=0$.
On peut m\^eme imposer que $c=2\pi (n+2)$.
\end{lemme}
On indique simplement comment se ramener \`a
des r\'esultats d'E.~Giroux.

D'apr\`es \cite{[Gi1]}, quitte \`a effectuer
des isotopies $C^0$-petites sur les tores $T^2 \times \lbrace k\rbrace$,
on peut toujours supposer que
le feuilletage $\xi T^2 \times \lbrace k \rbrace$
comporte un nombre fini d'orbites p\'eriodiques non d\'eg\'en\'er\'ees.
On se place dans cette situation.
La structure $\xi \vert_{T^2 \times [0,k]}$
est alors compl\`etement d\'ecrite (voir \cite{[Gi5]})
par ---~outre la torsion $n$ et  $\xi T^2 \times \lbrace i \rbrace$, pour $i=0,k$~---
la donn\'ee d'une famille d'anneaux $(A_i )_{1\leq i\leq k}$ deux \`a deux disjoints
dans $T^2 \times [0,k]$, s'appuyant sur les
orbites p\'eriodiques de $\xi T^2 \times \lbrace k \rbrace$.
On consid\`ere le produit $T^2 \times [0,k]$
comme naturellement inclus dans $T^2 \times [0,k+1]$.
Toujours d'apr\`es \cite{[Gi5]}, plonger
$(T^2 \times [0,k] ,\xi )$ dans
$(T^2 \times [0 ,2\pi (n+2)  ],
\ker (\cos (\theta +\theta_0 ) dx+\sin (\theta +\theta_0 )dy))$
revient a trouver une famille d'anneaux $(A_i' )_{1\leq i\leq l}$ deux \`a deux disjoints
 dans $T^2 \times [k,k+1]$
dont les bords co\"\i ncident avec ceux des anneaux de la
famille $(A_i )_{1\leq i\leq l}$ et
telle que la r\'eunion des deux familles constitue un ou deux
tores incompressibles dans $T^2 \times [0,k+1]$.

En d'autres termes, la preuve du lemme~\ref{lemme : normalisation}
se r\'eduit \`a la d\'emonstration du lemme suivant,
laiss\'ee \`a la sagacit\'e du lecteur~:

\begin{lemme} : Soit $(\alpha_i )_{1\leq i\leq l}$ une
famille d'arcs deux \`a deux disjoints inclus
dans $S^1 \times [0,1]$, dont les extr\'emit\'es sont
incluses
dans $S^1 \times \lbrace 1\rbrace$.
Il existe une famille d'arcs $(\alpha_i' )_{1\leq i\leq l}$
inclus dans $S^1 \times [1,2]$, dont les extr\'emit\'es
co\"\i ncident avec celles des arcs de la famille $(\alpha_i )_{1\leq i\leq l}$,
et telle que la r\'eunion des familles $(\alpha_i )_{1\leq i\leq l}$
et $(\alpha_i' )_{1\leq i\leq l}$ forme un ou deux
cercles non contractiles dans $S^1 \times [0,2]$.
\end{lemme}

$\;$\\
{\bf \'Elimination des singularit\'es}

\begin{lemme}\label{lemme : elimination} :
Soient $\xi$ une structure universellement
tendue sur $V=T^2 \times R$, et
$A=S^1 \times \lbrace * \rbrace \times [-1,1] \subset T^2 \times [-1,1]
\subset T^2 \times R$. On suppose que $\partial A$ est transverse \`a $\xi A$
avec un signe constant
et que $\xi A$ ne contient que des singularit\'es de m\^eme signe.
On supppose de plus que $ \xi T^2 \times \lbrace 1\rbrace$ est
un feuilletage sans singularit\'e ni composante de Reeb.

Il existe une isotopie $C^0$-petite de $T^2 \times [-1,1]$
fixant un voisinage de $A \cup T^2 \times \lbrace 1\rbrace$
telle que l'image de $T^2 \times \lbrace -1\rbrace$
porte un feuilletage sans singularit\'e ni composante de Reeb.
\end{lemme}
\begin{demonstration} : On utilise
des techniques d\'evelopp\'ees dans \cite{[Gi1]}.
On se place sur un rev\^etement
de degr\'e deux $p:V'\rightarrow V$ de $V$ dans lequel la pr\'eimage de
$A$ est constitu\'ee de deux anneaux disjoints $A_1$ et $A_2$.
En particulier, dans ce rev\^etement, $p^{-1} (A)$
d\'ecoupe $p^{-1} (T^2 \times [-1,1] )$ en deux tores
solides, dont un est d\'enot\'e $T$.
Le bord de $T$ est constitu\'e de la r\'eunion de $p^{-1} (A)$
et de deux anneau $B_{-1}$ et $B_1$ qui sont envoy\'es par $p$
respectivement sur
$T^2 \times \lbrace - 1 \rbrace$ et $T^2 \times \lbrace 1 \rbrace$.
On sait par hypoth\`ese que, si $\xi'$ d\'esigne le rappel de $\xi$
dans $V'$,
le feuilletage $\xi' (A_2 \cup B_1 )$ ne porte que des singularit\'es
de m\^eme signe et est soit sortant, soit rentrant le long du bord.
En particulier, $l(\partial (A_2 \cup B_1 ))=0$.
On a donc aussi $l(\partial (A_1 \cup B_{-1})) =0$.
Comme $\xi' A_1$ ne pr\'esente que des singularit\'es
de m\^eme signe et que $\xi' A_1$ est
soit rentrant, soit sortant le long de $\partial A_1$,
ces singularit\'es s'\'eliminent toutes par une isotopie de $A_1$ fixant
$B_{-1}$ (voir \cite{[Gi1]}).
On en d\'eduit que les singularit\'es de $B_{-1}$ s'\'eliminent
par une isotopie $C^0$-petite de $B_{-1}$ fixant un voisinage
de son bord (voir \cite{[Gi1]}), ce qui,
transcrit dans $V$, fournit le r\'esultat recherch\'e.
\end{demonstration}

\subsection{Structures universellement tendues sur $R\times S^1 \times [0,1]$}
On \'etend la notion de torsion aux structures
d\'efinies sur, respectivement, $R\times S^1 \times S^1 =\lbrace (x,y,\theta )\rbrace$
et $R\times S^1 \times [0,1]$ (avec des coordonn\'ees similaires),
en rempla\c cant la recherche de plongements de
$T^2 \times [0,2\pi ]$ par celle de plongements (propres et $\pi_1$-injectifs)
de $(R\times S^1 \times [0,2\pi] ,\ker (\cos n\theta dx
+\sin n\theta dy))$ respectivement non s\'eparant et
parall\`eles au bord. On parle alors de torsion annulaire.
Le r\'esultat suivant d\'ecoule
de \cite{[Co3]}.
Il peut \^etre \'egalement d\'eduit des
techniques d\'evelopp\'ees dans \cite{[Gi5]}.

\begin{lemme}\label{lemme : annulaire}{\rm \cite{[Co3],[Gi5]}} :
Toute structure de contact universellement tendue sur
$R\times S^1 \times S^1$ est de torsion annulaire finie.
\end{lemme}

\begin{corollaire}\label{lemme : fini} :
Toute structure de contact universellement tendue $\xi$
sur $V=R\times S^1 \times [-1,1]=\lbrace (x,y,\theta )\rbrace$ qui a pour
\'equation $\cos f(\theta )dx+\sin f(\theta )dy =0$ sur
$R\times S^1 \times [ -1 ,-1 +\epsilon ]$
et $R\times S^1 \times [ 1-\epsilon ,1]$
a une torsion annulaire finie.
\end{corollaire}
\begin{demonstration} :
On note $g:[2,3] \rightarrow R$ une fonction de
d\'eriv\'ee strictement positive avec $g(2)=f(1)\; mod\; 2\pi$
et $g(3)=f(-1)\; mod\; 2\pi$.
Le champ de plans $\zeta$ donn\'e par l'\'equation $\cos g(\theta )dx +\sin \theta dy=0$
est  une structure de contact sur $R\times S^1 \times [2,3]=\lbrace (x,y,\theta )\rbrace$.
En identifiant $R\times S^1 \times \lbrace 1\rbrace \subset V$
avec $R\times S^1 \times \lbrace 2\rbrace$ d'une part, et
$R\times S^1 \times \lbrace -1\rbrace \subset V$
avec $R\times S^1 \times \lbrace 3\rbrace \subset V$ d'autre part,
on obtient, en prolongeant $\xi$ par $\zeta$,
une structure de contact $\eta$ sur $W=R\times S^1 \times S^1$.
La structure $\eta$ est universellement tendue par application
d'une version annulaire du th\'eor\`eme de
recollement~\ref{theoreme : chirurgie} (voir \cite{[Co2]}).
On applique alors le lemme~\ref{lemme : annulaire} pour conclure.
\end{demonstration}

\subsection{Structures  tendues sur le
demi-espace $R^2 \times [0,+\infty[$}
 Sur le demi-espace $R^2 \times [0,+\infty [$, la notion
de torsion dispara\^\i t comme le montre le th\'eor\`eme
de classification suivant, d\^u \`a Y.~Eliashberg.

\begin{lemme}\label{lemme : demi-espace}{\rm \cite{[ElT]}} : Si $\xi$ est une structure
tendue sur $R^2 \times [0,+\infty [=\lbrace (x,y,t)\rbrace$
qui a pour \'equation $\cos tdx+\sin tdy=0$
sur $R^2 \times [0,\epsilon [$, alors
elle est conjugu\'ee \`a la structure d\'efinie globalement
par la m\^eme \'equation par un diff\'eomorphisme qui est
l'identit\'e pr\`es du bord.
\end{lemme}

\section{D\'ecomposition topologique du probl\`eme}

Une surface $S$ $(\neq S^2 ,D^2)$ plong\'ee dans une vari\'et\'e de
dimension $3$ est dite {\it incompressible}
si son groupe fondamental s'injecte dans celui de $V$.

Soit $V$ une vari\'et\'e de dimension $3$
orientable, irr\'eductible et close qui contient un tore incompressible.

On sait d'apr\`es W.~Jaco, P.~Shalen et K.~Johannson~\cite{[Ja], [JS], [Jo]} qu'il existe une
collection minimale finie $(T_i)_{0\leq i\leq n}$ de tores
incompressibles deux \`a deux disjoints,
unique \`a permutation et isotopie pr\`es, telle
que toute composante de $V\setminus \cup_{0\leq i\leq n} T_i$
soit ou un fibr\'e de Seifert, ou atoro\"\i dale (c'est-\`a-dire dans laquelle
tout tore incompressible est parall\`ele \`a une composante de bord).

Si $V$ est un fibr\'e de Seifert, on note $T$ un tore
incompressible quelconque de $V$.
Dans le cas contraire, la collection $(T_i)_{0\leq i\leq n}$
d\'ecoupant $V$ est non vide et on pose $T=T_0$.
On note $\pi :\bar{V} \rightarrow V$ le rev\^etement de
$V$ de groupe $\pi_1 (T)$.
Deux situations peuvent se produire :
\begin{itemize}
\item a) tous les relev\'es de $T$ dans $\bar{V}$ sauf un sont des plans~;
\item b) il existe un rev\^etement $p: \tilde{V} \rightarrow V$
de groupe $Z$ dans lequel deux relev\'es $\tilde{T_1}$ et $\tilde{T_2}$
de $T$ sont des anneaux conjugu\'es \`a $S^1 \times R$.
\end{itemize}
Dans les deux cas, on note $\bar{T}$ un relev\'e compact de
$T$ dans $\bar{V}$.

\begin{lemme}\label{lemme : topo1} :
Dans le cas a), on note $\bar{V_1}$ et $\bar{V_2}$
l'adh\'erence dans $\bar{V}$ des deux composantes
de $\pi^{-1} (V\setminus T)$ adjacentes \`a $\bar{T}$.
La vari\'et\'e $\bar{V}$ est alors obtenue en recollant
\`a chaque composante de $\partial (\bar{V_1} \cup \bar{V_2} )$
un demi-espace ($\simeq R^2 \times [0,\infty [$).
En particulier, $\bar{V}$ est diff\'eomorphe \`a l'int\'erieur
de $\bar{V_1} \cup \bar{V_2}$ dans $\bar{V}$.
Topologiquement, elle est \'egalement conjugu\'ee
\`a $T^2 \times R$.
\end{lemme}
\begin{demonstration} :
Comme $V$ est Haken,
elle est rev\^etue par $R^3$ et tout relev\'e $P$ de $T$
dans $R^3$ d\'ecoupe $R^3$ en deux demi-espaces.
En recollant, suivant les cas, un de ces deux demi-espaces
aux composantes de bord non compactes de
$\bar{V_1}$ et $\bar{V_2}$, on obtient
un rev\^etement de $V$, qui est de groupe
$\pi_1 (T)$ d'apr\`es le th\'eor\`eme de  Van Kampen.
Il est donc conjugu\'e \`a $\bar{V}$.

Le fait que $\bar{V}$ soit conjugu\'e \`a
$T^2 \times R$ est un r\'esultat classique de
topologie expliqu\'e dans \cite{[Si],[Ja]}.
\end{demonstration}
\vskip 0.5cm

De la m\^eme mani\`ere, on montre le r\'esultat suivant~:

\begin{lemme}\label{lemme : topo2} : Dans le cas b), la vari\'et\'e $\tilde{V}$
est conjugu\'ee \`a $R\times S^1 \times R$ et les deux anneaux $\tilde{T_1}$
et $\tilde{T_2}$
\`a $R\times S^1 \times \lbrace -1\rbrace$ et $R\times S^1 \times \lbrace 1\rbrace$.
On note $\tilde{V_1}$ et $\tilde{V_2}$  l'adh\'erence des composantes
de $p^{-1} (V\setminus T)$
adjacentes \`a, respectivement, $\tilde{T_1}$,
et $\tilde{T_2}$ et qui ne rencontrent pas le produit
$R\times S^1 \times ]-1,1[$.
Trois cas peuvent alors se produire~:
\begin{itemize}
\item $b_1 )$ La composante de Seifert adjacente \`a $T$
ne fibre pas au-dessus d'une bande de M\"oebius et, pour $i=1,2$, les composantes de
$\partial \tilde{V_i} \setminus \tilde{T_i}$
sont toutes des plans.
La vari\'et\'e $\tilde{V}$ est alors obtenue \`a partir de $R\times S^1 \times [-1,1]$
en recollant
$\tilde{V_1}$ sur $T^2 \times \lbrace - 1\rbrace$ et $\tilde{V_2}$ sur
$T^2 \times \lbrace  1\rbrace$, puis en recollant \`a chaque
composante de bord restante de $\tilde{V_1}$ et $\tilde{V_2}$
des demi-espaces ($\simeq R^2 \times [0,\infty[$).
Dans ce cas, on note $q : \tilde{V} \rightarrow \bar{V}$
l'application de rev\^etement qui
envoie $\tilde{T_1}$ sur $\bar{T}$.
L'image de $\tilde{V_1}$ par $q$ est conjugu\'ee \`a $T^2 \times [0,+\infty [$
et $q$ induit un diff\'eomorphisme de $\tilde{V_2}$ sur son image.
 Toutes les
composantes de bord de $q(\tilde{V_1} \cup R\times S^1 \times [-1,1] \cup \tilde{V_2} )$
sont des plans et on a~:
$\bar{V} \simeq Int(q(\tilde{V_1} \cup (R\times S^1 \times [-1,1]) \cup \tilde{V_2} ) )$.

\item $b_2 )$ La composante de Seifert
$M$ adjacente \`a $T$ fibre
au-dessus d'un ruban de M\"oebius
et deux relev\'es cons\'ecutifs $\bar{T}$ et $\bar{T'}$ de $T$ dans $\bar{V}$ sont des tores.
On note $\bar{V_1}$ et $\bar{V_2}$
les deux composantes de $\pi^{-1} (V\setminus T)$ attach\'ees
\`a l'ext\'erieur du produit $P$ relevant $M$ et
d\'elimit\'e dans $\bar{V}$ par ces tores.
Les composantes de bord de $\bar{V_1}$ et $\bar{V_2}$ autres que
$\bar{T}$ et $\bar{T}'$ sont toutes des plans.
La vari\'et\'e $\bar{V}$ est alors conjugu\'ee \`a
$Int(\bar{V_1} \cup P \cup \bar{V_2} )$.
\item $b_3 )$ Il existe un rev\^etement de $V$
conjugu\'e \`a $R\times S^1 \times S^1$ dans lequel
$T$ poss\`ede un rappel conjugu\'e \`a $R\times S^1 \times \lbrace *\rbrace$.
\end{itemize}
\end{lemme}
\begin{demonstration} :
Dire que deux anneaux proprement plong\'es et $\pi_1$-inject\'es
dans $R \times S^1 \times R$ bordent un produit est un r\'esultat
classique pour lequel on renvoie \`a \cite{[Wa]}.

Pour obtenir une image plus pr\'ecise de la situation, on
va utiliser la position de $T$ par rapport \`a
la d\'ecomposition de Jaco-Shalen de $V$.

Tout d'abord, on peut
dire que si $V$ est un fibr\'e de Seifert, on est dans le cas $b_3 )$.
Dans cette situation en effet, $V$ poss\`ede  un rev\^etement fini $\hat{V}$ qui fibre sur le cercle.
Quitte \`a consid\'erer un rev\^etement de degr\'e $2$ de $\hat{V}$,
on se ram\`ene au cas o\`u $\hat{V}$ fibre
au-dessus d'une surface orientable $S$.
Tout rappel $\hat{T}$ de $T$ dans $\hat{V}$ est alors
incompressible et est, \`a isotopie pr\`es,
satur\'e pour cette fibration (cf. \cite{[He]}). La surface $S$ est donc de genre
strictement positif.
Soit $\alpha$ une courbe ferm\'ee simple dans $S$ telle que $\hat{T}$
soit isotope au rappel de $\alpha$ dans $\hat{V}$.
On note $\beta$ une courbe ferm\'ee simple dans $S$
que l'on ne peut pas disjoindre de $\alpha$ par isotopie,
dont l'existence est assur\'ee par la minoration du genre de $S$,
et $\hat{T}'$ le rappel de $\beta$ dans $\hat{V}$.
Il n'existe alors pas d'isotopie de disjonction
de $\hat{T}$ et $\hat{T}'$ (et $\hat{T}$
est un tore normal).
Le rev\^etement $W$ de $\hat{V}$ de groupe $\pi_1 (\hat{T}' )$, vu
comme rev\^etement de $V$, poss\`ede les propri\'et\'es
requises dans le cas $b_3 )$ (voir par exemple \cite{[Co3]}).

Sinon, on est par hypoth\`eses dans le cas o\`u $T$ est isotope
\`a un des tores de la d\'ecomposition de Jaco-Shalen.

Soit $\bar{T}'$ un relev\'e de $T$ dans
$\bar{V}$, distinct de $\bar{T}$, qui est un anneau infini (conjugu\'e \`a
$R\times S^1 $) ou un tore. On note $K$ l'adh\'erence de la composante de
$\bar{V} \setminus (\bar{T} \cup \bar{T}' )$ qui rencontre
$\bar{T}$ et $\bar{T}'$.
D'apr\`es le th\'eore\`eme de Van Kampen, $\pi_1 (K) =\pi_1 (\bar{V} )=Z^2$.
Il existe alors dans $\bar{T}'$ une courbe ferm\'ee simple homologue
\`a une courbe ferm\'ee simple de $\bar{T}$.
Ces deux courbes bordent une surface dans $K$
que l'on peut toujours rendre incompressible
et qui est donc un anneau, car $\pi_1 (K )=Z^2$.
 D'apr\`es le th\'eor\`eme de l'anneau~VIII.10 de \cite{[Ja]},
$\bar{T}$ et $\bar{T'}$ doivent border un m\^eme
relev\'e d'une composante du d\'ecoupage de Jaco-Shalen
qui est un fibr\'e de Seifert $M$.
En particulier, tous les relev\'es de $T$ dans $\bar{V}$
qui ne sont pas des plans doivent border l'adh\'erence
d'une m\^eme composante de $\pi^{-1} (V\setminus T)$.\\
\\
$\bullet$ Si $\bar{T}'$ est un tore, la composante de Seifert $M$
poss\`ede un rev\^etement de degr\'e deux conjugu\'e
\`a $T^2 \times [0,1]$~: il s'agit d'une fibration
en cercles sur un ruban de M\"oebius (cf. \cite{[Ja]}).
Tous les relev\'es de $T$ autres que $\bar{T}$ et $\bar{T}'$
sont, d'apr\`es ce qui pr\'ec\`ede, des plans.
La m\^eme \'etude que celle de la configuration $a)$
permet de terminer la d\'emonstration dans ce cas.\\
\\
$\bullet$ Si $\bar{T'}$ est un anneau, par passage au rev\^etement, on obtient
que tous les relev\'es de $T$ dans $\tilde{V}$
qui ne sont pas des plans doivent border l'adh\'erence
d'une m\^eme composante de $p^{-1} (V\setminus T)$, ce qui suffit
\`a ramener la d\'emonstration, comme dans le cas du lemme pr\'ec\'edent,
\`a une utilisation du th\'eor\`eme de Van Kampen et \`a des raisonnements
\'el\'ementaires sur les rev\^etements : par exemple,
en recollant des demi-espaces, judicieusement
choisis dans le rev\^etement universel de $V$, aux composantes
de bord de $\tilde{V_1} \cup (R\times S^1 \times [-2,2]) \cup \tilde{V_2}$
(qui sont toutes planes d'apr\`es ce qui pr\'ec\`ede),
on obtient un rev\^etement de $V$ de groupe $Z$ qui est donc conjugu\'e \`a $\tilde{V}$.
\end{demonstration}
\begin{remarque}{\rm : Plus g\'en\'eralement, on peut dire que,
si
$T$ est normal (ce qui \'equivaut \`a dire que $T$ n'est pas isotope
\`a l'un des tores de la d\'ecomposition de Jaco-Shalen-Johannson), on est dans le cas $b_3 )$~:
il existe en effet un relev\'e annulaire $\tilde{T_3}$ de $T$ situ\'e hors du
produit d\'elimit\'e dans $\tilde{V}$ par $\tilde{T_1}$ et
$\tilde{T_2}$ (voir \cite{[Co3]}), et on est dans la situation ou au moins
deux anneaux parmi ceux de la famille $(\tilde{T_i} )_{1\leq i\leq 3}$
bordent un produit $P$ (conjugu\'e \`a $R\times S^1 \times [0,1]$)
avec une coorientation (obtenue en relevant une coorientation de $T$)
rentrante dans le produit pour l'un et sortante pour
l'autre. Si $\phi$ est un diff\'eomorphisme entre ces deux anneaux
commutant avec $p$ et pr\'eservant l'orientation, la vari\'et\'e obtenue
par identification des deux bords de $P$ \`a l'aide
de $\phi$ fournit une vari\'et\'e $W$, conjugu\'ee
\`a $R\times S^1 \times S^1$, qui est un rev\^etement de
$V$ dans lequel $T$ poss\`ede un rappel (donn\'e par
l'image de $\partial P$) conjugu\'e \`a $R\times S^1 \times \lbrace *\rbrace$.
La d\'emonstration du th\'eor\`eme~\ref{theoreme}
dans ce cas a d\'ej\`a \'et\'e
essentiellement trait\'ee dans \cite{[Co3]}.}
\end{remarque}

\section{Construction de la suite d'exemples}

Pour construire des structures de contact tendues,
on dispose d'un important travail effectu\'e
par D.~Gabai~\cite{[Ga]} d'une part et par Eliashberg-Gromov~\cite{[El3],[Gr]}
et Eliashberg-Thurston~\cite{[ElT]}
d'autre part.

Un feuilletage est dit {\it tendu} s'il poss\`ede une transversale
ferm\'ee qui rencontre toutes ses feuilles.
L'existence des feuilletages tendus a \'et\'e largement
\'etudi\'ee par D.~Gabai.
On rappelle qu'une surface incompressible
est dite {\it minimale} si elle minimise
le genre dans sa classe d'homologie.
Ici, la surface n'est pas suppos\'ee connexe,
mais on suppose qu'aucune composante n'est une
sph\`ere ou un disque, et le genre d\'esigne alors
la somme des genres des composantes connexes.

\begin{theoreme}\label{theoreme : Gabai}{\rm \cite{[Ga]}} :
Soient $V$ une vari\'et\'e compacte irr\'eductible de dimension $3$
bord\'ee par une r\'eunion non vide de tores et $S$ une
surface minimale dans $V$ qui repr\'esente un \'el\'ement
non nul de $H_2 (V,\partial V,Z)$.
Il existe un feuilletage tendu dont une feuille est $S$
et qui trace sur $\partial V$ un feuilletage
sans singularit\'e ni composante de Reeb.
\end{theoreme}

Un {\it feuilletact} est un champ de plans $\xi =\ker \alpha$
tel que la $3$-forme $\alpha \wedge d\alpha$ soit
de signe constant. Cette notion interpole entre
les notions de feuilletage (o\`u $\alpha \wedge d\alpha$ est
identiquement nulle) et de structure de contact (o\`u
$\alpha \wedge d\alpha$ ne s'annulle pas).

\begin{theoreme}\label{theoreme : limite}{\rm \cite{[ElT]}} : Tout feuilletact $\xi$ est limite $C^0$
de structures de contact qui co\"\i ncident avec $\xi$
sur un compact o\`u il est d\'ej\`a de contact.
\end{theoreme}

Pour prouver qu'une vari\'et\'e de contact $(V,\xi )$
est  tendue, on dispose
d'un crit\`ere d\^u \`a Y.~Eliashberg et M.~Gromov~\cite{[El3],[Gr]}~:
il suffit qu'il existe une vari\'et\'e symplectique $(W,\tilde{\omega} )$
($\tilde{\omega}$ est une $2$-forme ferm\'ee non d\'eg\'en\'er\'ee)
qui {\it remplisse} $(V,\xi )$, ce qui signifie que
$\partial W=V$, que $\tilde{\omega} \vert_{\xi}$ est
non d\'eg\'en\'er\'ee et que $\tilde{\omega}$ oriente $V=\partial W$
comme $\xi$.

En particulier (voir \cite{[ElT]}), il suffit qu'il existe une
$2$-forme ferm\'ee $\omega$ sur $V$
qui {\it domine} $\xi$, c'est-\`a-dire
telle que $\omega \vert_\xi$ soit non d\'eg\'en\'er\'ee,
ainsi qu'une structure n\'egative $\xi'$ \'egalement
domin\'ee par $\omega$.

Dans ce cas en effet, la $2$-forme $\tilde{\omega}=p^* \omega +\epsilon
ds\alpha$ obtenue
\`a partir de $\omega$ sur $V\times [0,1]$ ---~$s$ est la coordonn\'ee
sur $[0,1]$, $\epsilon$ est un r\'eel positif
assez petit, $p: V\times [0,1] \rightarrow V$ la projection
sur $V$ et $\alpha =0$ est une \'equation de $\xi$~--- d\'etermine un remplissage symplectique
de $(V,\xi ) \coprod (V,\xi' )$
et le th\'eor\`eme de remplissage de \cite{[El3],[Gr]} assure alors que $\xi$
et $\xi'$ sont tendues.

Plus g\'en\'eralement, dans le cas o\`u
$(V, \xi )=\partial (W, \tilde{\omega} )$ n'est pas compacte,
il faut supposer de plus,
outre le fait que $\tilde{\omega} \vert_{\xi} >0$,
qu'il existe une structure presque complexe $J$
sur $W$ qui pr\'eserve $\xi$ et telle que $\tilde{\omega} (*,J*)$
soit une m\'etrique riemanienne $g$ avec les propri\'et\'es suivantes~:
\begin{itemize}
\item $g$ est compl\`ete~;
\item le rayon d'injectivit\'e de $g$ est minor\'e par un r\'eel strictement
positif~;
\item la courbure sectionnelle de $g$ est major\'ee.
\end{itemize}
On dit alors que $(W,\tilde{\omega} )$ \`a une {\it g\'eom\'etrie finie \`a l'infini},
et, d'apr\`es \cite{[El3],[ElT],[Gr]},
toute vari\'et\'e de contact bord\'ee par une vari\'et\'e
symplectique qui poss\`ede une g\'eom\'etrie finie \`a
l'infini est tendue.

D'apr\`es un travail de D.~Sullivan~\cite{[Su]},
imposer l'existence d'une telle $2$-forme pour un feuilletage
\'equivaut \`a dire qu'il est tendu.

Une autre strat\'egie pour construire des structures
tendues est de d\'eterminer leur
comportement par rapport \`a certaines op\'erations de chirurgies.

\begin{theoreme}\label{theoreme : chirurgie}{\rm \cite{[Co2]}} :
Soient $(V,\xi )$ une vari\'et\'e de contact
et $T\subset V$ un tore incompressible.
Si $\xi T$ est un feuilletage lin\'eaire et
si la vari\'et\'e $(V\setminus T, \xi \vert_{V\setminus T} )$
est universellement tendue, alors $(V,\xi )$ l'est aussi.
La conclusion est la m\^eme si $T$ est un anneau incompressible
plong\'e qui poss\`ede un voisnage tubulaire
conjugu\'e \`a~:
$(R\times S^1 \times [-\epsilon ,\epsilon ]=\lbrace (x,y,\theta )\rbrace ,
\ker (\cos (\theta +\theta_0 )+\sin (\theta +\theta_0 )))$, avec $\epsilon$,
$\theta_0 \in R$.
\end{theoreme}

C'est sur cette cha\^\i ne de r\'esultats
qu'on se fonde pour prouver le lemme suivant,
l\'eg\`ere adaptation d'un \'enonc\'e
de K.~Honda, W.~Kazez et G.~Mati\'c \cite{[HKM]}, qui
pr\'ecise le travail d'Eliashberg et Thurston.

\begin{lemme}\label{lemme : minimal}{\rm (voir \cite{[HKM]})} :
Soit $V$ une vari\'et\'e orient\'ee de dimension
$3$, qui est  irr\'eductible, compacte et dont le bord est
une r\'eunion non vide de tores incompressibles.
Soit $S$ une surface minimale plong\'ee dans $V$,
chaque composante connexe de $S$ \'etant de bord non vide plong\'e
dans $\partial V$.
Il existe une structure de contact universellement tendue positive
$\xi$ sur $V$, qui trace sur $\partial V$ un feuilletage
lin\'eaire et sur chaque
composante de $S$ un feuilletage transverse au bord avec un signe constant, sans orbite
p\'eriodique et dont toutes les singularit\'es sont de m\^eme
signe.
\end{lemme}
\begin{demonstration} : D'apr\`es le th\'eor\`eme~\ref{theoreme : Gabai}
de Gabai, il existe
un feuilletage tendu $F$ dont $S$ est une feuille
et qui trace sur $\partial V$ un feuilletage sans
singularit\'e ni composante de Reeb.
Quitte \`a modifier $F$, on peut supposer qu'il est
conjugu\'e a un feuilletage produit $S\times [-1,1]$
pr\`es de $S\simeq S\times \lbrace 0\rbrace$.

Soit alors $\omega$ une $2$-forme ferm\'ee,
donn\'ee par \cite{[Su]}, qui domine $F$.
On modifie $\omega$ pr\`es de chaque composante $T$ de $\partial V$
de la mani\`ere suivante.
Soit $T\times [-1,0] =\lbrace (x,y,\theta )\rbrace$
un syst\`eme de coordonn\'ees
pr\`es de $T\simeq T\times \lbrace 0\rbrace$
tel que $F$ soit tangent \`a
$\partial_\theta$ et transverse
\`a $\partial_y$ ($F\cap T$ est sans composante de Reeb).
On suppose de plus que $S\times [-1,1] \cap T$ est donn\'e
par l'\'equation $dy=0$.
La $2$-forme $\omega$ est ferm\'ee et s'\'ecrit donc dans
ces coordonn\'ees~: $\omega =a dx\wedge dy +d\eta$
pour un certain r\'eel $a$ et une certaine $1$-forme
$\eta$.

Si $g :[-1,0] \rightarrow R$ est une fonction
positive, nulle pr\`es de $-1$ et strictement
positive sur $[-\frac12 ,0]$, qui vaut $1$ au voisinage de $0$,
et si $\chi :[-1,0] \rightarrow R$
vaut $1$ sur $[-1,-\frac14 ]$ et $0$ pr\`es de $0$,
on pose, pour $t\in R$, $$\omega_t =adx\wedge dy + d(\chi \eta )
+tg(\theta )d\theta \wedge dx.$$
Cette forme $\omega_t$ est ferm\'ee, co\"\i ncide avec
$\omega$ pr\`es de $T\times \lbrace -1\rbrace$
et, si $t$ est assez grand, domine $F$.
On fixe un tel r\'eel $t$, assez grand pour chaque composante de bord,
 et on note $\omega_t'$
la $2$-forme obtenue en raccordant $\omega_t$ \`a $\omega$
pr\`es de chaque composante de bord.
On colle alors \`a chaque composante  de bord $T$
un produit $T\times [0,+\infty [=\lbrace (x,y,\theta )\rbrace$
(les coordonn\'ees prolongent celles
d\'ej\`a construites sur $T\times [-1,0]$)
sur lequel
on prolonge $F$ par un feuilletage produit
et $\omega_t'$ par la $2$-forme $adx\wedge dy +td\theta\wedge dx$.
On note $V'$ cette nouvelle vari\'et\'e, $F'$ le
nouveau feuilletage et $\omega_t"$ la
nouvelle forme.

On revient \`a pr\'esent sur $V$.
On trace sur $S$ un feuilletage
orient\'e $L$ transverse \`a $\partial S$,
dont les singularit\'es sont des selles
ou des foyers de divergence positive (pour une certaine
orientation de $S$ fix\'ee) et
ne pr\'esentant pas d'orbite p\'eriodique.
On choisit $L$ sortant le long de $\partial S$.
Un tel feuilletage est donn\'e par le noyau d'une
$1$-forme $\beta$ avec $d\beta >0$
(voir par exemple \cite{[Gi2]}).
Soit alors $f : [-1, 1] \rightarrow R$ une fonction lisse,
nulle en $\pm 1$ et strictement positive \`a l'int\'erieur.
Si $t$ d\'esigne la coordonn\'ee transverse \`a $S$
dans son voisinage $S\times [-1,1]$, la $1$-forme
$\alpha =dt+f(t) \beta$ v\'erifie~: $\alpha  \wedge d\alpha
= f(t) dt\wedge d\beta$.
Le noyau de $\alpha$ est donc une structure de contact
sur $S\times ]-1,1[$ qui trace $L$ sur $S$, et co\"\i ncide
avec $F$ au bord. On construit ainsi sur $V$ un feuilletact
positif $\xi_0$ qui trace le feuilletage d\'esir\'e sur $S$.
La proximit\'e de $\xi_0$ \`a $F$ est donn\'ee par la taille de $f$.

Soit \`a nouveau $T$ une composante de $\partial V$
et $\lbrace x=x_0 \rbrace$ une \'equation de
$S\cap T$. Le feuilletage $\xi_0 T$
a une \'equation de la forme $dy+h(x,y)dx=0$,
avec $h(x_0 ,y) >0$.
Soit alors $H : T^2 \times [0,\infty [ \rightarrow R$
une fonction lisse qui poss\`ede les propri\'et\'es suivantes~:
\begin{itemize}
\item $\frac{\partial H}{\partial \theta} (x,y,\theta )>0$~;
\item $\lim_{\theta \rightarrow +\infty} H=c_0$, $c_0 \in R$~;
\item $H(x,y,1)=c_1$, $c_1 \in R$~;
\item $H(x,y,0)=h(x,y)$.
\end{itemize}
Le champ de plans d\'efini sur $T^2 \times [0,+\infty [$
comme $\lbrace dy+Hdx=0\rbrace$ est une structure de contact
qui prolonge $\xi_0$ en un feuilletact.
Elle est transverse \`a $\partial_y$
et trace sur chaque tore $T\times \lbrace 1\rbrace$
un feuilletage lin\'eaire d'\'equation $dy+c_1 dx=0$.
Elle trace de plus sur $(S\cap T)\times [0,1]$
un feuilletage non singulier,
dirig\'e par $\partial_\theta$, et donc sans orbite p\'eriodique.

On effectue ce prolongement pr\`es de chaque composante
de $\partial V$ pour obtenir un feuilletact sur $V'$.
Si $t$ est choisi assez grand, la $2$-forme $\omega_t"$
domine ce champ de plans.
Le feuilletact ainsi construit sur $V'$ peut \^etre
approxim\'e, d'apr\`es le th\'eor\`eme~\ref{theoreme : limite},
par une structure de contact
$\xi$.
Comme la condition de domination est ouverte
pour la topologie $C^0$, la structure $\xi$ est domin\'ee par $\omega_t"$.
Pour construire une structure de contact n\'egative $\xi'$ domin\'ee par
$\omega_t"$, il suffit d'appliquer \`a nouveau le th\'eor\`eme~\ref{theoreme : limite}
pour approximer le feuilletage $F'$ par une structure n\'egative.
La vari\'et\'e $V'$ n'est pas compacte,
mais la forme $\omega_t"$ est ``constante" \`a l'infini.
On peut appliquer le th\'eor\`eme de Gromov et Eliashberg
\cite{[Gr],[El3]}
pour conclure que $(V',\xi )$ est  universellement tendue~: on s'aper\c coit
que la forme symplectique
$\tilde{\omega}=p^* \omega"_t +\epsilon ds\alpha$ obtenue
\`a partir de $\omega"_t$ sur $V\times [0,1]$ ---~$s$ est la coordonn\'ee
sur $[0,1]$, $\epsilon$ est un r\'eel positif
assez petit, $p: V\times [0,1] \rightarrow V$ la projection
sur $V$ et $\alpha =0$ est une \'equation de $\xi$~---, ainsi
que son rappel dans le rev\^etement universel de $V\times [0,1]$,
a une g\'eom\'etrie
finie \`a l'infini (cf. \cite{[ElT]}).

Pour conclure, il suffit de remarquer que $V$ est conjugu\'ee
\`a $V'\setminus (\cup_{T\subset \partial V} T\times ]1,\infty [)$.
 \end{demonstration}

\vskip 0.5cm

Soit $V$ une vari\'et\'e irr\'eductible, orient\'ee, close,
$T\subset V$ un tore incompressible  et $T \times [-1,1]$,
$T\times \lbrace 0\rbrace \simeq T$ un voisinage tubulaire de $T$.
On note $S\subset V\setminus T \times ]-1,1[$
une surface minimale (dont l'existence est assur\'ee
dans \cite{[He]}) qui rencontre les deux tores $T \times \lbrace \pm 1\rbrace$
et dont toutes les composantes connexes sont \`a bord.

D'apr\`es le lemme~\ref{lemme : minimal}, il existe une structure
de contact universellement tendue
$\xi^0$ sur $V\setminus T \times ]-1,1[$
avec les propri\'et\'es suivantes~:
\begin{itemize}
\item $\xi^0$ trace un feuilletage lin\'eaire
sur $T \times \lbrace \pm 1\rbrace$~;
\item $\xi^0$ est transverse \`a $\partial S$
et trace sur $S$ un feuilletage sans orbite p\'eriodique dont
toutes les singularit\'es ont le m\^eme signe.
\end{itemize}
Sur $T \times [-1,1]$, muni des coordonn\'ees $(x,y,\theta )$,
on prolonge $\xi^0$ par une structure
d'\'equation $\cos f_n (\theta )dx+\sin f_n (\theta )dy=0$,
o\`u $f_n [-1,1] \rightarrow R$ est une fonction de d\'eriv\'ee
strictement positive qui v\'erifie~:
$f_n (-1 )\in [0, 2\pi [$~; $f_n (1)\in [2n\pi ,2(n+1) \pi [$,
ces valeurs \'etant telles que les feuilletages
trac\'es sur $T \times \lbrace \pm 1\rbrace$
par les structures situ\'ees de part et d'autre co\"\i ncident.

D'apr\`es le th\'eor\`eme de recollement~\ref{theoreme : chirurgie},
toute structure $\xi_n$ ainsi construite sur $V$ est universellement tendue.

\begin{theoreme}\label{theoreme : np} : La suite $(\xi_n )_{n\in N}$
comporte une infinit\'e de structures de contact
deux \`a deux non isomorphes.
\end{theoreme}

On note $C$ la classe d'isotopie du tore $T$ dans $V$.
Pour montrer le th\'eor\`eme~\ref{theoreme : np}, il
suffit de d\'emontrer le r\'esultat suivant~:

\begin{theoreme}\label{theoreme : intermediaire} : $Tor (V,\xi_n ,C)<\infty$.
\end{theoreme}

Voici, en effet, comment d\'eduire le th\'eor\`eme~\ref{theoreme : np}
du th\'eor\`eme~\ref{theoreme : intermediaire}.

D'apr\`es \cite{[Co3]}, si $\xi$ est une structure universellement tendue sur
$V$ seul un nombre fini (\`a reparam\'etrisation pr\`es) de classes d'isotopie $D$ de
plongements de $T^2 \times [0,2\pi ]$ dans $V$ fournissent une torsion $Tor (V,\xi ,D)$
non nulle.
On note $(C^n_i )_{1\leq i\leq k_n}$ ces classes pour $(V,\xi_n )$, $C=C_0^n$, et
$T_n =\sup_{1\leq i\leq k_n} \lbrace Tor (V,\xi_n ,C^n_i ),\; avec\;  Tor (V,\xi_n ,C^n_i )<\infty \rbrace$,
qui est bien d\'efini et sup\'erieur \`a $n$ car $n\leq Tor (V,\xi_n ,C^n_0 )<\infty$.
Ainsi, si $p\in N$ est choisi assez grand,
$T_p >T_n$.
Si $u_n$ est une suite d'entiers telle
que la suite $T_{u_n}$ soit strictement croissante,
on est alors assur\'e que les structures
de contact $\xi_{u_n}$ sont deux \`a deux non isomorphes.

\section{Preuve du th\'eor\`eme~\ref{theoreme : intermediaire}}
\subsection{Configuration a)}
Dans cette partie, on va montrer que $Tor (V,\xi_n ,C)$
vaut $n-1$, $n$, $n+1$ ou $n+2$.

Soient $T\subset V$ un tore incompressible dans la configuration a)
et $\xi_n$ une structure de contact universellement tendue
construite selon la m\'ethode propos\'ee plus
haut. On note $\bar{\xi_n}$ le rappel de $\xi_n$ sur
$\bar{V}$ et $\bar{T} \times [-1,1]$ le relev\'e de $T\times [-1,1]$
qui contient $\bar{T} \simeq \bar{T} \times \lbrace 0\rbrace$.
On note de plus $\bar{V_1}'$ et $\bar{V_2}'$ les deux relev\'es
de $V\setminus T\times ]-1,1[$ dans $\bar{V}$ adjacents
\`a $\bar{T} \times [-1,1]$.

\begin{lemme}\label{lemme : reduction} : $(\bar{V} ,\bar{\xi_n} )$ est contactomorphe
\`a $(Int (\bar{V_1}' \cup (\bar{T} \times [-1,1]) \cup \bar{V_2}'  ) ,
\bar{\xi_n} \vert_{Int(\bar{V_1}' \cup (\bar{T} \times [-1,1]) \cup \bar{V_2}' ) })$.
\end{lemme}
\begin{demonstration} : On sait par construction
que les relev\'es non compacts de $T\times \lbrace \pm 1\rbrace$
dans $\bar{V}$ ont des voisinages tubulaires disjoints
conjugu\'es \`a
$(R^2 \times [-\alpha,\alpha ], \ker (\cos tdx+\sin tdy))$.

On remarque alors \`a l'aide du lemme~\ref{lemme : demi-espace}
(et du lemme~\ref{lemme : topo1})
que la structure $\bar{\xi_n}$ est conjugu\'ee, sur chaque composante
de $\bar{V} \setminus Int (\bar{V_1}' \cup (\bar{T} \times [-1,1]) \cup \bar{V_2}')$, \`a
$(R^2 \times [0 ,\infty [,\ker ( \cos tdx+\sin tdy))$.
Dans $\bar{V}$, on a donc un plongement de
$(R^2 \times [-\alpha ,\infty [, \ker (\cos tdx+\sin tdy))$
pr\`es de chaque composante de bord non compacte de $\bar{V_1}'$
et $\bar{V_2}'$ (celles-ci \'etant
identifi\'ees avec $R^2 \times \lbrace 0\rbrace$).
On conclut  gr\^ace \`a un changement
de variable explicite bien connu
(voir par exemple~\cite{[Co2]}, remarque~4.3), ou par une nouvelle
application du lemme~\ref{lemme : demi-espace}, qui
envoie $(R^2 \times [-\alpha ,\infty [, \ker (\cos tdx+\sin tdy))$
sur $(R^2 \times [-\alpha ,0[, \ker (\cos tdx+\sin tdy))$
en laissant invariant un voisinage de $R^2 \times \lbrace -\alpha \rbrace$.
\end{demonstration}

\hskip 0.5cm

Par la suite, $M$ d\'esigne la vari\'et\'e $\bar{V_1}'$ ou $\bar{V_2}'$
priv\'ee de ses composantes de bord non compactes, et
$\zeta$ la structure induite par
$\bar{\xi_n}$ sur $M$. On rappelle que $M$ est diff\'eomorphe
\`a $T^2 \times [0,+\infty [$.

\begin {lemme}\label{lemme : morceau} : La vari\'et\'e
$(M,\zeta )$ est de torsion nulle.
\end{lemme}
\begin{demonstration} : Quitte \`a changer
l'orientation de $S$, on suppose
que toutes les singularit\'es de $\zeta S$
sont positives.
On note $S'$ un rappel de $S$ dans $M$, de groupe
fondamental $Z$. La surface $S'$ est diff\'eomorphe
\`a un demi-anneau  infini $S^1 \times [0,\infty [$.
Son bord est inclus dans celui de $M$.
Le feuilletage de $S'$ est form\'e
 de singularit\'es positives, et toute feuille
de ce feuilletage aboutit en un temps fini \`a l'une des singularit\'es.
Il contient un unique cycle ferm\'e, parall\`ele au bord.

Par une isotopie $C^0$-petite de $S'$
on fait appara\^\i tre, sur chaque liaison
entre deux selles de $\zeta S'$, un foyer et
une selle positives.
On note $S"$ cette nouvelle surface.
Dans $\zeta S"$ on ne peut pas trouver de
chemin form\'e de liaisons contenant $4$ selles cons\'ecutives.

Pour tout anneau $A\subset S"$, $\partial S"\subset \partial A$,
l'ensemble $\omega_-$-limite de $A$ par le flot de $\zeta A$
est un graphe de singularit\'es  fini $\Gamma$, form\'e d'arbres
venant se greffer sur le cycle ferm\'e.
Soit $s$  une selle de ce graphe et $l$ une s\'eparatrice
stable de $s$. La s\'eparatrice $l$ provient alors soit
d'un foyer, soit d'une selle $s'$. \`A leur tour,
les s\'eparatrices stables de $s'$ proviennent
soit d'un foyer, soit d'une selle dont
chaque s\'eparatrice stable provient alors d'un foyer.

Le graphe $\Gamma$ est donc inclus dans un graphe $\Gamma'$ fini
dont toutes les extr\'emit\'es sont des foyers
et qui contient toutes les s\'eparatrices
stables de ses sommets.
Le graphe $\Gamma'$ poss\`ede alors un voisinage tubulaire
(qui est un anneau) de bord transverse \`a $\zeta A$.

 On d\'eduit  de ces remarques
qu'il existe une exhaustion de $S"$ par une
suite d'anneaux embo\^\i t\'es
$(A_i )_{i\in N}$, telle que~:
\begin{itemize}
\item pour tout $i\in N$, $\partial S"\subset \partial A_i$~;
\item pour tout $i\in N$, $\partial A_i$ est transverse \`a $\zeta$
et le feuilletage $\zeta A_i$ est sortant sur les deux bords.
\end{itemize}
Comme $A$ est proprement plong\'e et $\pi_1$-inject\'e dans $M$,
on peut montrer, \`a l'aide des techniques d\'evelopp\'ees dans \cite{[Wa]},
que le couple $(M,A)$ est conjugu\'e
\`a $(S^1 \times S^1 \times [0,+\infty [,\lbrace *\rbrace \times S^1 \times [0,+\infty [)$.
En particulier,
il existe une exhaustion de $M$ par une famille de tores
\'epais (conjugu\'es \`a $T^2 \times [0,1]$) $(T_i )_{i\in N}$
telle que~:
\begin{itemize}
\item $\partial M\subset \partial T_i$~;
\item $T_i$ rencontre transversalement $S"$ le long de $A_i$.
\end{itemize}
D'apr\`es le lemme~\ref{lemme : elimination} d\'eduit
des travaux de E.~Giroux, on peut supposer de plus que
$\partial T_i$ pr\'esente un feuilletage
sans singularit\'e ni composante
de Reeb.

Pour conclure, on remarque que chaque anneau $A_i$ porte un arc legendrien
$\gamma_i$, d\'efini comme r\'eunion de feuilles de $\zeta A_i$,
 qui joint un bord de $A_i$ \`a l'autre, et
tel que $tb(\gamma_i ,A_i ) = -\frac12$ car toutes les singularit\'es de $\zeta A_i$
ont le m\^eme signe.
En particulier, le lemme~\ref{lemme : 0} permet de conclure
que la torsion de $T_i$ est nulle pour tout $i$.
C'est donc aussi le cas pour $M$.
\end{demonstration}

\vskip 0.5cm

On applique le lemme~\ref{lemme : addition} d'addition des torsions
\`a la d\'ecomposition~:
$$Int(\bar{V_1}' \cup (\bar{T} \times [-1,1]) \cup \bar{V_2}' ).$$
On note que la torsion de $(\bar{T} \times [-1,1] ,\bar{\xi_n} \vert_{\bar{T} \times [-1,1]} )$
est $n-1$ ou $n$ d'apr\`es le lemme~\ref{lemme : torsion} et que celle des deux autres composantes
de la d\'ecomposition est nulle d'apr\`es le lemme~\ref{lemme : morceau}.
On en d\'eduit que la torsion de
$Int(\bar{V_1}' \cup (\bar{T} \times [-1,1]) \cup \bar{V_2}' ),
\bar{\xi_n} \vert_{Int(\bar{V_1}' \cup (\bar{T} \times [-1,1]) \cup \bar{V_2}' ) })$
vaut $n-1$, $n$, $n+1$ ou $n+2$.
D'apr\`es le lemme~\ref{lemme : reduction}, c'est aussi
le cas pour $(\bar{V} ,\bar{\xi_n} )$.

On conclut en remarquant que~:
$$n-1\leq Tor (V,\xi_n , C)\leq Tor(\bar{V} ,\bar{\xi_n} ).$$

\subsection{Configuration b)}
\subsubsection{Cas $b_1 )$}

Cette fois, on suppose qu'il existe un rev\^etement $p:\tilde{V} \rightarrow V$
de $V$ conjugu\'e au produit $R\times S^1 \times R$ o\`u $T$ poss\`ede
deux relev\'es  $\tilde{T_1}$ et $\tilde{T_2}$
qui sont des anneaux.
Dans la construction de $\xi_n$ sur $V$, on part d'un
voisinage tubulaire $T\times [-1,1]$ de $T$.
Dans  $\tilde{V}$, on note
$\tilde{T_1} \times [-1,1 ]$
et $\tilde{T_2} \times [-1,1 ]$ les deux voisinages
de $\tilde{T_1}$ et $\tilde{T_2}$ qui rel\`event $T\times [-1,1]$.
D'apr\`es le lemme~\ref{lemme : topo2},
on peut supposer que $\tilde{V}=R\times S^1 \times R$ avec, pour tout
$t\in [-1,1]$,
$\tilde{T_1}\times \lbrace t\rbrace =R\times S^1 \times \lbrace -2 +t\rbrace$ et
$\tilde{T_2}\times \lbrace t\rbrace =R\times S^1 \times \lbrace 2+t\rbrace$.
On note $\tilde{V_1}'$ et $\tilde{V_2}'$ les deux composantes
de $p^{-1} (V\setminus (T\times ]-1,1[))$ adjacentes \`a, respectivement,
$\tilde{T_1} \times [-1,1]$ et $\tilde{T_2} \times [-1,1]$ et qui ne rencontrent pas le produit
$R\times S^1 \times ]-3,3[$.

On remarque que pour $i=1,2$, il existe $j_i \in\lbrace 1,2\rbrace$
tel que $\tilde{V_i}'$ rev\^ete $\bar{V_{j_i}}'$.
On note $p_{j_i}$ l'application de rev\^etement.

Comme, hormis  $\tilde{T_i} \times \lbrace \pm 1\rbrace$, les composantes de
$\partial \tilde{V_i}'$ sont des plans, la m\^eme d\'emonstration que celle du
lemme~\ref{lemme : reduction}, bas\'ee sur le lemme~\ref{lemme : demi-espace} de classification
des structures tendues sur le demi-espace, fournit que les sous-vari\'et\'es
$(R\times S^1 \times ]-\infty ,-3] ,\tilde{\xi_n} )$
et $(R\times S^1 \times [3,+\infty [ ,\tilde{\xi_n} )$ de
$\tilde{V}=R\times S^1 \times R$ sont
conjugu\'ees aux vari\'et\'es
$(\tilde{V_1}' ,\tilde{\xi_n} \vert_{\tilde{V_1}'} )$ et
$(\tilde{V_2}' ,\tilde{\xi_n} \vert_{\tilde{V_2}'} )$ priv\'ees de leurs
composantes de bord non annulaires,
et donc  que~: $$(\tilde{V} ,\tilde{\xi_n} )\simeq (Int (\tilde{V_1}' \cup
(R\times S^1 \times [-3,3]) \cup
\tilde{V_2}' ),\tilde{\xi_n} \vert_{Int (\tilde{V_1}' \cup (R\times S^1 \times [-3,3]) \cup
\tilde{V_2}' )} ).$$
On note $q$, comme dans le lemme~\ref{lemme : topo2},
 l'application de rev\^etement de $\tilde{V}$ sur $\bar{V}$ qui
envoie
$\tilde{T_1}$ sur $\bar{T}$.
De la m\^eme fa\c con que pr\'ec\'edemment, comme,
d'apr\`es le lemme~\ref{lemme : topo2}, toutes les composantes de $\partial
(q(Int (\tilde{V_1}' \cup
(R\times S^1 \times [-3,3]) \cup
\tilde{V_2}' ))$ sont des plans, on obtient que~:
$$(\bar{V} ,\bar{\xi_n} )\simeq (q(Int (\tilde{V_1}' \cup
(R\times S^1 \times [-3,3]) \cup
\tilde{V_2}' )),q_* \tilde{\xi_n} \vert_{Int (\tilde{V_1}' \cup (R\times S^1 \times [-3,3]) \cup
\tilde{V_2} )}' ).$$

\begin{remarque}{\rm : On prendra garde que, si $q(\tilde{V_1}' )=\bar{V_1}'$
ou $\bar{V_2}'$,
on n'a pas en revanche $q(\tilde{V_2}' )=\bar{V_2}'$ ou $\bar{V_1}'$.}
\end{remarque}

Pour simplifier les notations, on
suppose que $q(\tilde{V_1}' )=\bar{V_1}'$, et donc en particulier
que $j_1 =1$.

On veut montrer que la torsion est finie sur $\bar{V}$.
D'apr\`es l'isomorphisme pr\'ec\'edent, il suffit de
consid\'erer un produit $P=T^2 \times [0,2\pi ]$ plong\'e
de mani\`ere $\pi_1$-injective dans $(q(Int (\tilde{V_1}' \cup
(R\times S^1 \times [-3,3]) \cup
\tilde{V_2}' )),q_* \tilde{\xi_n} \vert_{Int (\tilde{V_1}' \cup (R\times S^1 \times [-3,3]) \cup
\tilde{V_2}' )} )$.

Soient $\bar{V_1}"$ et $\bar{V_2}"$, les sous-vari\'et\'es de
$\bar{V}$ obtenues en retirant \`a $\bar{V_1}'$ et
$\bar{V_2}'$ les composantes de bord autres que $\bar{T}\times \lbrace \pm 1\rbrace$.
On rappelle que, d'apr\`es \cite{[Si]},
pour $i=1,2$, $\bar{V_i}" \simeq T^2 \times [0,+\infty[$.
D'apr\`es l'\'etude du cas a), et
notamment le lemme~\ref{lemme : morceau}, la torsion est
nulle sur $\bar{V_i}"$.

On remarque alors que si $k_1 ,k_2 \in N$ sont assez grand, $P\cap q(\tilde {V_i}' )$
est inclus dans $T^2 \times [0,k_1 ] \subset \bar{V_1}"\simeq T^2 \times [0,+\infty [$
si $i=1$, et dans $p_{j_2}^{-1} (T^2\times [0,k_2] )\subset q(\tilde{V_2} )\simeq \tilde{V_2}$ pour
$i=2$, o\`u $T^2\times [0,k_2] \subset \bar{V_{i_2}}" \simeq T^2 \times [0,+\infty [$.

De plus, d'apr\`es le lemme~\ref{lemme : normalisation},
il existe $(c_1 ,c_2 ,\theta_1 ,\theta_2 )\in R^4$ (ind\'ependant de $k_1$ et $k_2$)
tel que, pour $i=1,2$, $(T^2 \times [0,k_i ],\bar{\xi_n} \vert_{T^2 \times [0,k_i ]} )$
se plonge dans $(T^2 \times [0,c_i ],\ker (\cos (\theta +\theta_i ) dx+\sin (\theta +\theta_i )dy))$.

Ainsi, par passage au rev\^etement,
on en d\'eduit que
le rappel par $q$ de $P$ dans $\tilde{V}$
se plonge dans le mod\`ele
$$(W,\zeta )\simeq (R\times S^1 \times [0,c_1 ],
\ker (\cos (\theta +\theta_1 )dx+\sin (\theta +\theta_1 )  dy))
\cup (R \times S^1 \times [-3,3] ,\tilde{\xi_n} \vert_{R\times S^1 \times [-3,3]} )
$$ $$\cup (R\times S^1 \times [0,c_2 ] ,
\ker (\cos (\theta +\theta_2 )dx+\sin (\theta +\theta_2 ) dy))$$
(o\`u on identifie $R\times S^1 \times \lbrace 0 \rbrace$ avec
$R\times S^1 \times \lbrace -3 \rbrace$ d'une part, et
avec
$R\times S^1 \times \lbrace 3 \rbrace$ d'autre part).
La vari\'et\'e $(W,\zeta )$ est obtenue par collage
de trois vari\'et\'es de contact universellement
tendues le long d'anneaux ($R\times S^1 \times \lbrace \pm 3\rbrace$)
dans un voisinage tubulaire desquels la structure a par
construction une \'equation
du type $\cos (\theta +\theta_0 )dx +\sin (\theta +\theta_0 ) dy =0$,
o\`u $\theta \in [-\epsilon' ,\epsilon' ]$
d\'esigne la coordonn\'ee transverse.
En particulier, la version annulaire du th\'eor\`eme de
recollement~\ref{theoreme : chirurgie}
fournit que $(W,\zeta )$ est universellement tendue.

On est de plus bien dans
les hypoth\`eses du lemme~\ref{lemme : fini}~: sa conclusion
s'applique et la torsion annulaire
est finie sur $(W,\zeta )$.

On conclut la d\'emonstration du th\'eor\`eme~\ref{theoreme : intermediaire}.
en remarquant  que~:
$$Tor (V,\xi_n ,C)\leq Tor(\bar{V} ,\bar{\xi_n} )\leq Tor (W ,\zeta ).$$

\subsubsection{Cas $b_2 )$}
L'\'etude dans ce cas est essentiellement un m\'elange de celle
des configuration $a)$ et $b_1 )$, c'est pourquoi on se contente
d'en r\'esumer les \'etapes.
On note $\bar{V_1}'$ et $\bar{V_2}'$ les deux relev\'es
de $V\setminus (T\times ]-1,1[ )$ dans $\bar{V}$
inclus dans $\bar{V_1}$ et $\bar{V_2}$.
On note $P'$ le produit bord\'e dans $\bar{V}$
par les composantes toriques de $\partial \bar{V_i}'$.
Comme les composantes de bord de
$\bar{V_1}' \cup P' \cup \bar{V_2}'$
sont des plans, on obtient, comme dans l'\'etude du cas $a)$, un isomorphisme~:
$$(\bar{V} ,\bar{\xi_n} )\simeq (Int (\bar{V_1}' \cup P' \cup \bar{V_2}' ),\bar{\xi_n}
\vert_{ Int (\bar{V_1}' \cup P' \cup \bar{V_2}' )} ).$$

Pour conclure il faut rappeler que, pour $i=1,2$, la vari\'et\'e
obtenue en retirant \`a $\bar{V_i}'$ ses composantes de bord non compactes
est conjugu\'ee \`a $T^2 \times [0,+\infty [$ et
que la torsion y est finie (lemme~\ref{lemme : morceau}).
La torsion est \'egalement finie sur $P'$, par exemple par application du corollaire~\ref{lemme : fini}
\`a un rev\^etement cyclique de $P'$.
On applique alors le lemme d'addition
des torsions~\ref{lemme : addition} \`a la d\'ecomposition
ci-dessus pour obtenir que la torsion de $(\bar{V} ,\bar{\xi_n } )$
est finie, et donc aussi celle de $(V,\xi )$ dans la classe $C$ de $T$.

\subsubsection{Cas $b_3 )$}
Dans le cas $b_3 )$, on sait qu'il existe un rev\^etement $W$ de
$V$, conjugu\'e \`a $R\times S^1 \times S^1$,
dans lequel $T$ poss\`ede un relev\'e conjugu\'e
\`a $R\times S^1 \times \lbrace *\rbrace$.
Dans ce cas, le lemme~\ref{lemme : annulaire}
affirme que la torsion annulaire de $W$ est finie, ce qui
implique \`a nouveau imm\'ediatement le
th\'eor\`eme~\ref{theoreme : intermediaire}.

\vskip 1cm
\hskip 7cm Vincent Colin

\hskip 6cm Universit\'e de Nantes

\hskip 6cm D\'epartement de math\'ematiques

\hskip 6cm UMR 6629 du CNRS

\hskip 6cm 2, rue de la Houssini\`ere

\hskip 6cm BP 92208, 44322 Nantes cedex 3

\hskip 6cm Vincent.Colin@math.univ-nantes.fr
\end{document}